\RequirePackage{ifpdf}
\ifpdf
 	\documentclass[pdftex]{amsart}
 	\else
 	\documentclass[dvips]{amsart}
\fi

\usepackage{amsfonts,amsthm,latexsym,amsmath,amssymb,amscd,amsmath, mathrsfs, epsf}
\usepackage{graphicx}

\graphicspath{{./figures/}}
\ifpdf
	\usepackage{epstopdf}		
	\DeclareGraphicsExtensions{.png,.jpg,.eps,.epsf}
\fi

\newtheorem{theorem}{Theorem}
\newtheorem{proposition}[theorem]{Proposition}

\newtheorem{definition}{Definition}

\newcommand{\xvec}{\mathbf{x}}

\begin{document}
\title[Discriminants, symmetrized graph monomials, and SOS]{Discriminants, symmetrized graph monomials, and sums of squares certificates} 

\author[P.~Alexandersson]{Per Alexandersson}

\address{Department of Mathematics,
   Stockholm University,
   S-10691, Stockholm, Sweden}
\email{per@math.su.se}

\begin{abstract}
Here we present certificates for 5 classes of 6-edged multigraphs whose symmetrized graph monomials may be represented as sum of squares,
but not as linear combinations of partition square graphs. This is a complement to the results presented in \cite{sos}.
\end{abstract}

\maketitle

\section{Background}
In  \cite{sos}, we have the following definition:

\begin{definition}\label{def:graphmonomial}
Let $g$ be a directed graph, with vertices $x_1,\dots,x_n$ and adjacency matrix $(a_{ij}),$ being non-negative integers.
Define first its graph monomial $P_g$ as follows 
$$P_g(x_1,\dots,x_n):=\prod_{1\leq i , j \leq n} (x_i-x_j)^{a_{ij}},$$
where $a_{ij}$ is the number of directed edges joining $x_i$ with $x_j.$

The symmetrized graph monomial of $g$ is defined as
$$\tilde{g}(\xvec)=\sum_{\sigma \in S_n} P_g(\sigma \xvec),\quad \xvec = x_1,\dots,x_n.$$
\end{definition}
Two multigraphs are \emph{equivalent} if the first corresponding symmetrized graph monomial is 
a constant multiple of the other.

We study the set of symmetrized graph monomials obtained from multigraphs with 6 edges.
There are 212 such multigraphs.

One main result in \cite{sos} is the following proposition:
\begin{proposition}\label{pr:main5}
\rm {(i)} $102$ graphs with $6$ edges  have identically vanishing symmetrized graph monomial.  
\rm {(ii)} The remaining $110$ graphs are divided into 27 equivalence classes. 
\rm{(iii)} $12$ of these classes can be expressed as non-negative linear combinations of square graphs, 
i.e. lie in the convex cone spanned by the square graphs. 
\rm{(iv)} Of the remaining 15 classes, symmetrized graph monomial of $7$ of them change sign. 
\rm {(v)} Of the remaining 8 classes (which are presented in Fig.~\ref{fig:similarsixgraphs}) 
 the first 5 are sums of squares.
\end{proposition}

Proposition \rm{(i)},\rm{(ii), \rm{(iii)} and \rm{(iv)} follows from direct computation.

\section{Certificates}

Here we give certificates that the symmetrized graph monomial for the first $5$ classes given in Fig.~\ref{fig:similarsixgraphs} are sum of squares,
corresponding to the last case in Prop.~\ref{pr:main5}.
The graphs in each row yield the same polynomial, up to a constant. 
The symmetrized graph monomial from row $i$ is a constant multiple of the polynomial $v_i Q_i v_i^T,$
where $v_i$ is the coefficient vector and $Q_i$ is the corresponding symmetric positive semi-definite matrix, given below.
It is well-known that this certifies that the first 5 classes of graphs are sum of squares.

It is relatively straightforward to verify that the polynomials indeed are non-negative, 
using methods similarly to \cite{lax}.
By using Matlab together with Yalmip, one may verify that the last three classes cannot be expressed as sums of squares.

\begin{align*}
v_1=\{&x_1^2 x_2, x_1 x_2^2, x_1^2 x_3, x_2^2 x_3,x_1 x_3^2,x_2 x_3^2,x_1^2 x_4,x_2^2 x_4,x_3^2 x_4,x_1 x_4^2,x_2 x_4^2,x_3 x_4^2,\\
     &x_1^2 x_5,x_2^2x_5,x_3^2 x_5,x_4^2 x_5,x_1 x_5^2,x_2 x_5^2,x_3 x_5^2,x_4 x_5^2,x_5 x_6^2,x_5^2 x_6,x_4 x_6^2,x_4^2 x_6,\\
     &x_3 x_6^2,x_3^2 x_6,x_2 x_6^2,x_2^2x_6,x_1 x_6^2,x_1^2 x_6,x_1 x_2 x_3,x_1 x_2 x_4,x_1 x_3 x_4,x_2 x_3 x_4,x_1 x_2 x_5, \\
     &x_1 x_3 x_5,x_2 x_3x_5,x_1 x_4 x_5,x_2 x_4 x_5,x_3 x_4x_5,x_4 x_5 x_6,x_3 x_5 x_6,x_3 x_4 x_6,x_2 x_5 x_6,x_2 x_4 x_6, x_2 x_3 x_6, \\
     &x_1 x_5 x_6,x_1 x_4 x_6,x_1 x_3 x_6,x_1 x_2 x_6\}
\end{align*}

\begin{align*}
v_{2}=\{&x_5 x_6^2,x_5^2 x_6,x_4 x_6^2,x_4 x_5 x_6,x_4 x_5^2,x_4^2 x_6,x_4^2 x_5,x_3 x_6^2,x_3 x_5 x_6,x_3 x_5^2,x_3 x_4 x_6,\\
          &x_3 x_4 x_5,x_3x_4^2,x_3^2 x_6,x_3^2 x_5,x_3^2 x_4,x_2 x_6^2,x_2 x_5 x_6,x_2 x_5^2,x_2 x_4 x_6,x_2 x_4 x_5,x_2 x_4^2,\\
          &x_2 x_3 x_6,x_2 x_3 x_5,x_2 x_3 x_4,x_2x_3^2,x_2^2 x_6,x_2^2 x_5,x_2^2 x_4,x_2^2 x_3,x_1 x_6^2,x_1 x_5 x_6,x_1 x_5^2,\\
          &x_1 x_4 x_6,x_1 x_4 x_5,x_1 x_4^2,x_1 x_3 x_6,x_1 x_3 x_5,x_1x_3 x_4,x_1 x_3^2,x_1 x_2 x_6,x_1 x_2 x_5,x_1 x_2 x_4,\\
          &x_1 x_2 x_3,x_1 x_2^2,x_1^2 x_6,x_1^2 x_5,x_1^2 x_4,x_1^2 x_3,x_1^2 x_2\}
\end{align*}

\begin{align*}
v_{3}=\{&x_3 x_4^2,x_3^2 x_4,x_2 x_4^2,x_2 x_3 x_4,x_2 x_3^2,x_2^2 x_4,x_2^2 x_3,x_1 x_4^2,x_1 x_3 x_4,x_1 x_3^2,x_1 x_2 x_4,\\
          &x_1 x_2 x_3,x_1x_2^2,x_1^2 x_4,x_1^2 x_3,x_1^2 x_2\}
\end{align*}

\begin{align*}
v_{4}=\{&x_4^3,x_3 x_4^2,x_3^2 x_4,x_3^3,x_2 x_4^2,x_2 x_3 x_4,x_2 x_3^2,x_2^2 x_4,x_2^2 x_3,x_2^3,x_1 x_4^2,x_1 x_3 x_4,x_1 x_3^2,\\
          &x_1 x_2 x_4,x_1x_2 x_3,x_1 x_2^2,x_1^2 x_4,x_1^2 x_3,x_1^2 x_2,x_1^3\}
\end{align*}

\begin{align*}
v_{5}=\{&x_5 x_6^2,x_5^2 x_6,x_4 x_6^2,x_4 x_5 x_6,x_4 x_5^2,x_4^2 x_6,x_4^2 x_5,x_3 x_6^2,x_3 x_5 x_6,x_3 x_5^2,x_3 x_4 x_6,\\
           &x_3 x_4 x_5,x_3x_4^2,x_3^2 x_6,x_3^2 x_5,x_3^2 x_4,x_2 x_6^2,x_2 x_5 x_6,x_2 x_5^2,x_2 x_4 x_6,x_2 x_4 x_5,x_2 x_4^2,\\
           &x_2 x_3 x_6,x_2 x_3 x_5,x_2 x_3 x_4,x_2x_3^2,x_2^2 x_6,x_2^2 x_5,x_2^2 x_4,x_2^2 x_3,x_1 x_6^2,x_1 x_5 x_6,x_1 x_5^2,\\
           &x_1 x_4 x_6,x_1 x_4 x_5,x_1 x_4^2,x_1 x_3 x_6,x_1 x_3 x_5,x_1x_3 x_4,x_1 x_3^2,x_1 x_2 x_6,x_1 x_2 x_5,x_1 x_2 x_4,\\
           &x_1 x_2 x_3,x_1 x_2^2,x_1^2 x_6,x_1^2 x_5,x_1^2 x_4,x_1^2 x_3,x_1^2 x_2\}
\end{align*}

\begin{figure}[ht!]
\centering
\includegraphics[width=0.7\textwidth]{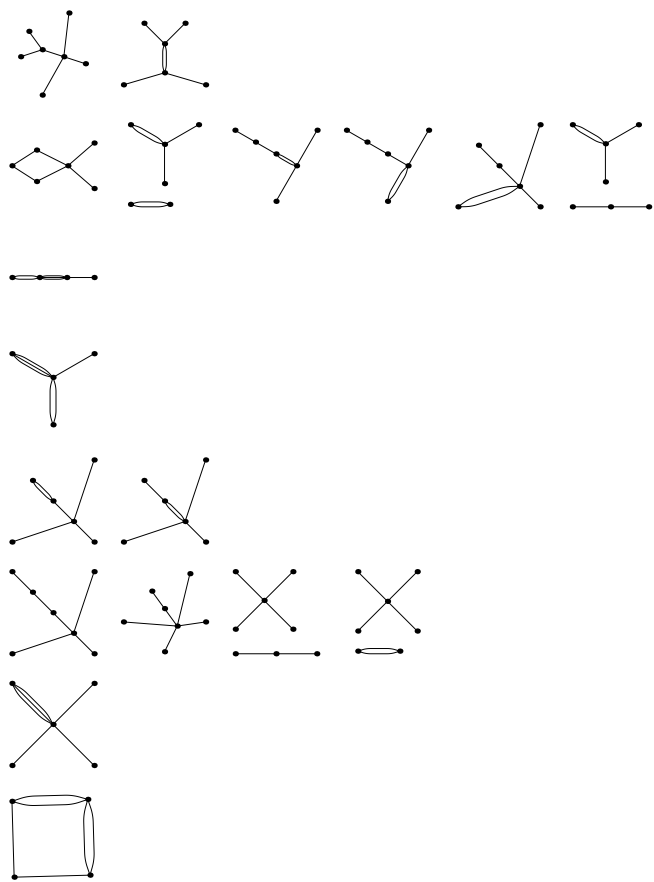}
\caption{$8$ equivalence classes of all non-negative graphs with $6$ edges.}\label{fig:similarsixgraphs}
\end{figure}

\begin{figure}[ht!]
\centering
\includegraphics[width=1.2\textwidth]{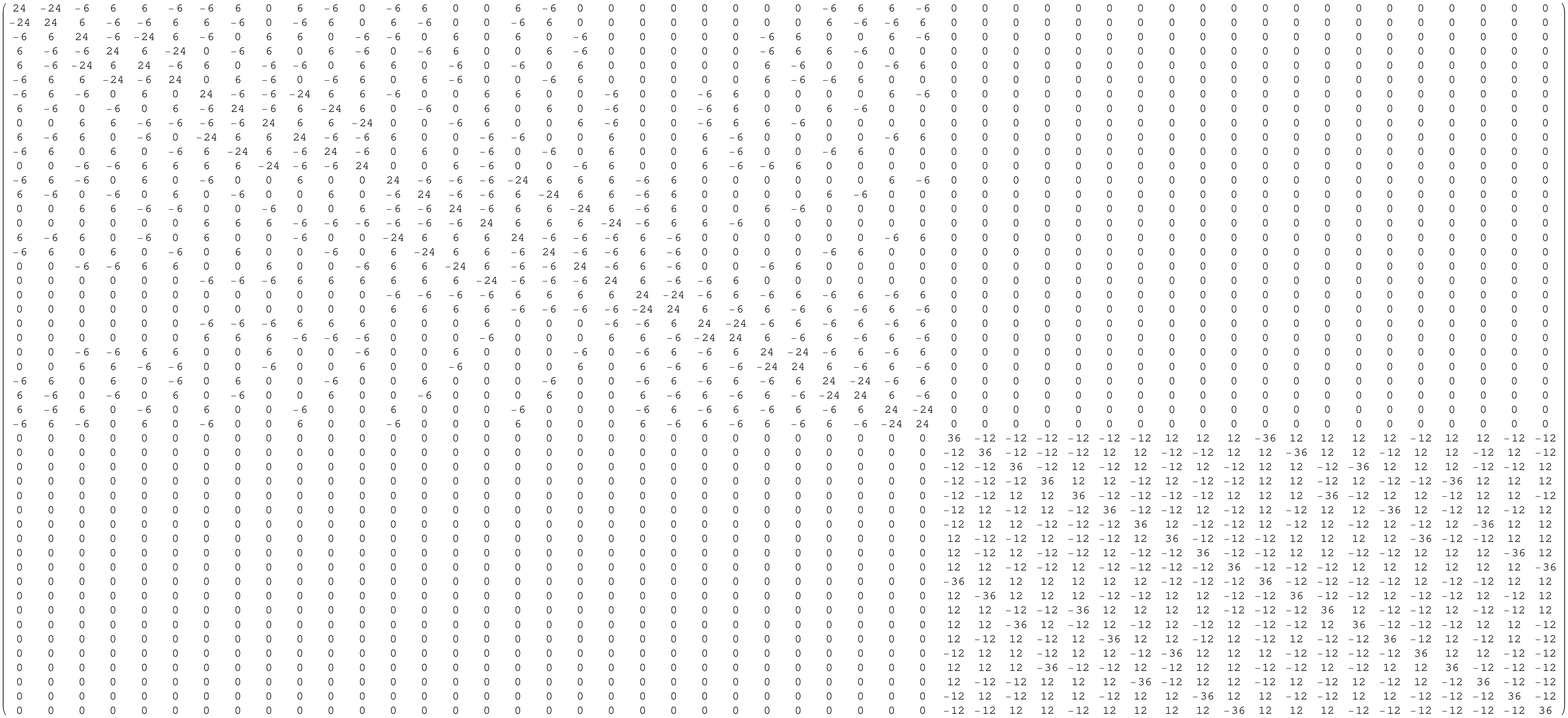}
\caption{$Q_1$}
\end{figure}

\begin{figure}[ht!]
\centering
\includegraphics[width=1.2\textwidth]{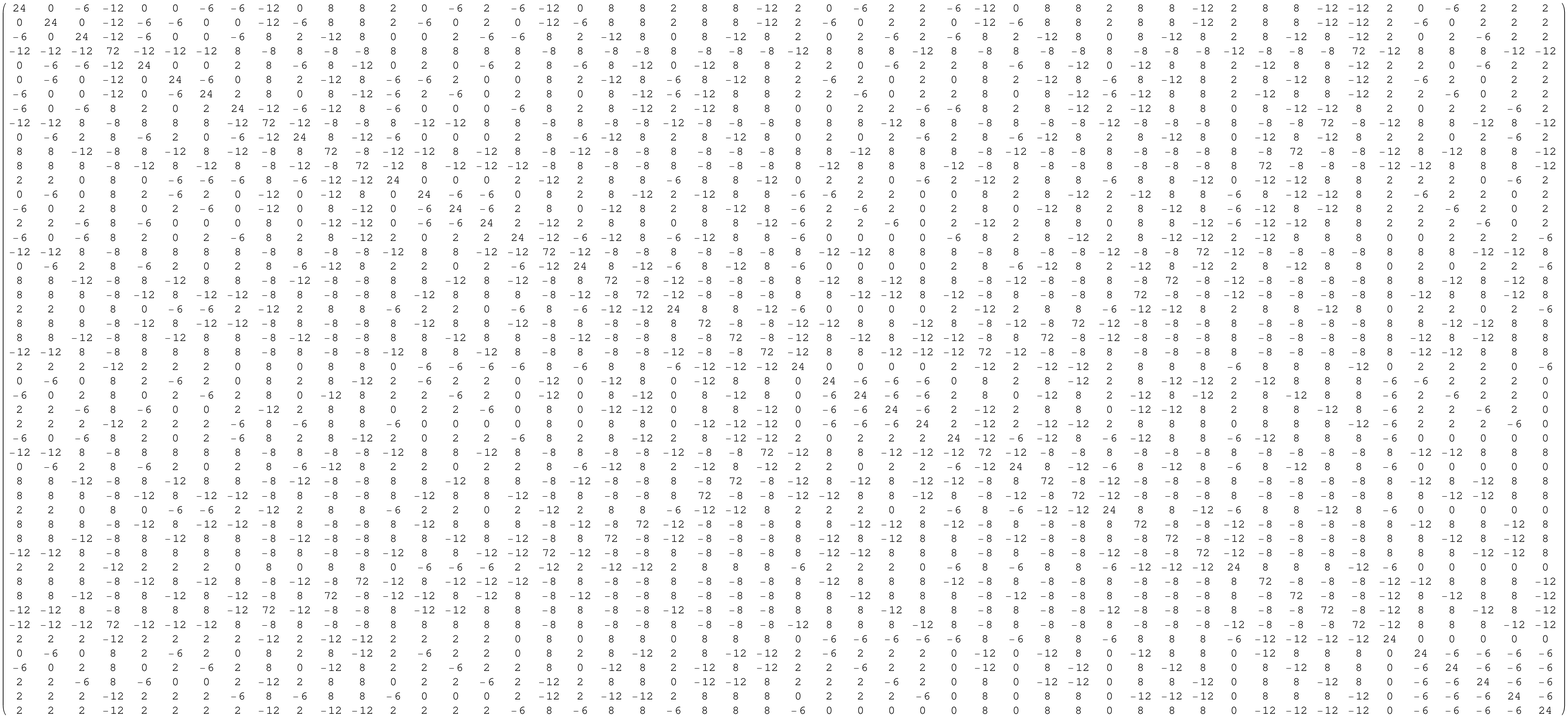}
\caption{$Q_{2}$}
\end{figure}

\begin{figure}[ht!]
\centering
\includegraphics[scale=1]{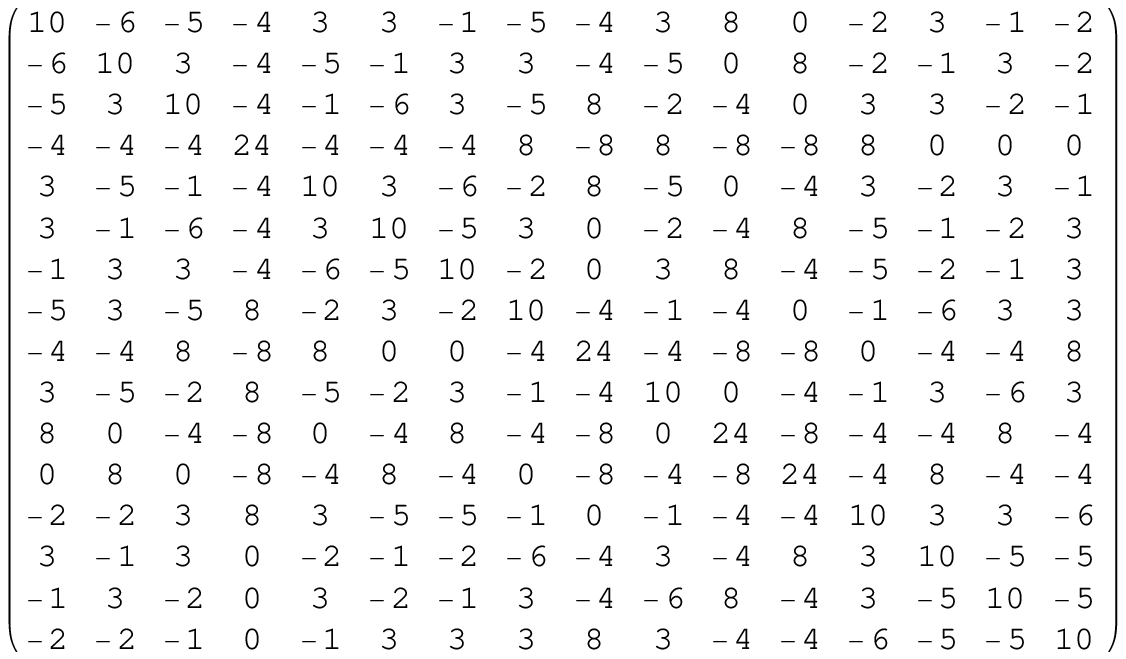}
\caption{$Q_{3}$}
\end{figure}

\begin{figure}[ht!]
\centering
\includegraphics[width=1.2\textwidth]{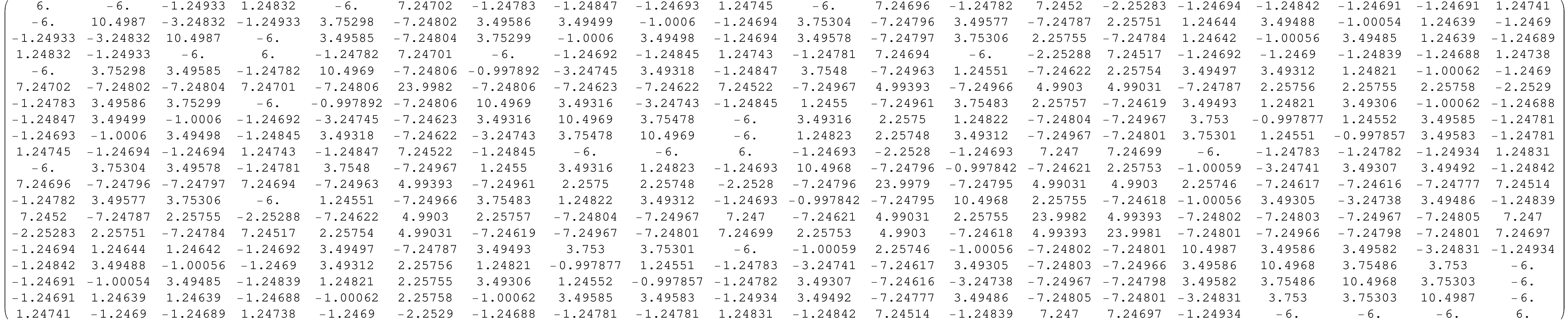}
\caption{$Q_{4}$}
\end{figure}

\begin{figure}[ht!]
\centering
\includegraphics[width=1.2\textwidth]{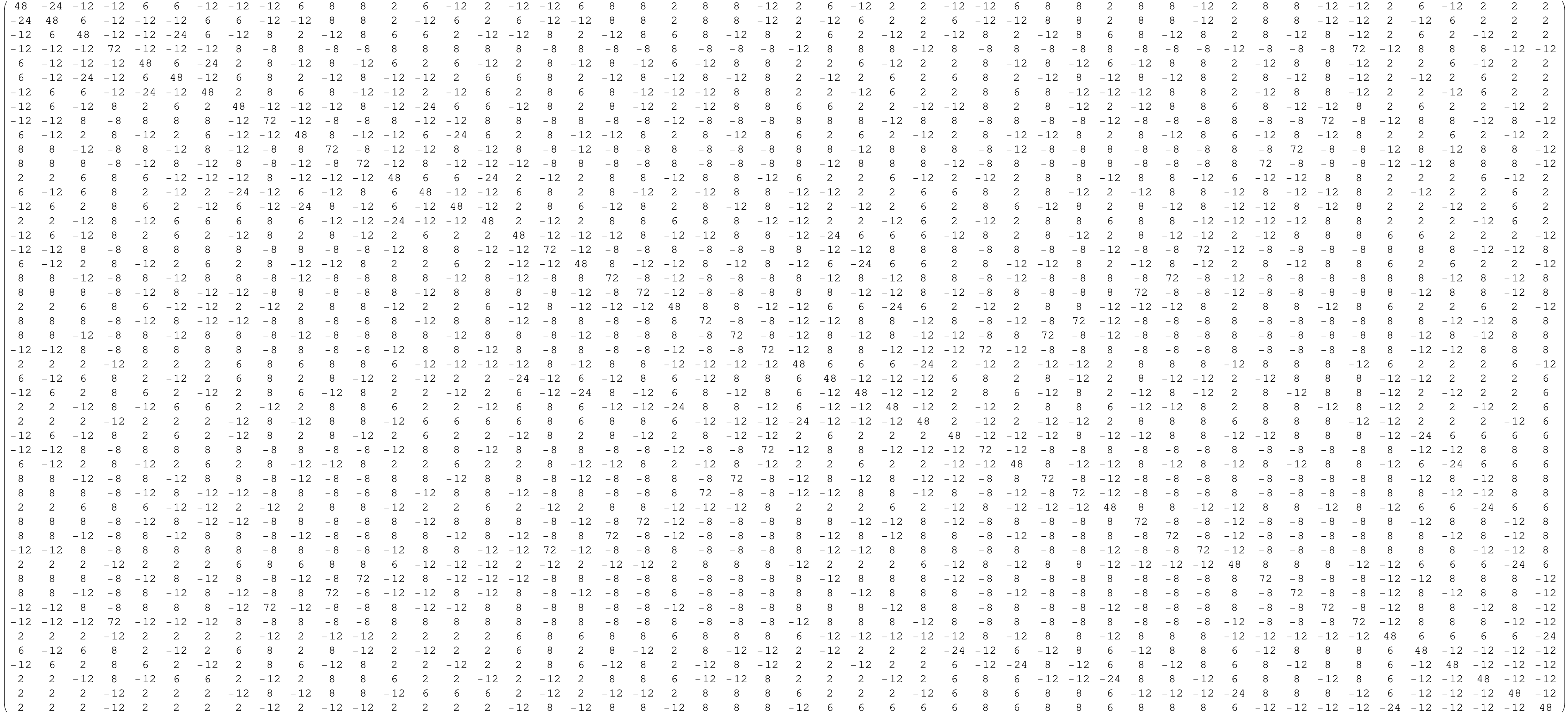}
\caption{$Q_{5}$}
\end{figure}

\end{document}